\documentclass[a4paper,12pt]{article}

\input{mymacros.sty}

\renewcommand{\bsa}{a}
\renewcommand{\bsb}{b}
\renewcommand{\bsc}{c}
\renewcommand{\bsd}{d}

\title{Dimer models and
homological mirror symmetry for triangles}
\author{Masahiro Futaki and Kazushi Ueda}
\date{}
\begin{document}

\maketitle

\begin{abstract}
We prove a conjecture on the relation
between dimer models, coamoebas and vanishing cycles
for the mirrors of two-dimensional toric Fano stacks
of Picard number one.
As a corollary, we obtain a torus-equivariant version
of homological mirror symmetry
for such stacks.
\end{abstract}

\section{Introduction}
With a convex lattice polygon $\Delta \subset \bR^2$
containing the origin in its interior,
one can associate a directed $A_\infty$-category
in three different ways:

\begin{itemize}
 \item
Let
$$
 W(x, y) = \sum_{(i, j) \in \bZ^2} a_{ij} x^i y^j
$$
be a Laurent polynomial
whose Newton polygon coincides with $\Delta$;
$$
 \Delta = \Conv \{ (i, j) \in \bZ^2 \mid a_{ij} \ne 0 \}.
$$
If the coefficients $a_{ij}$ are sufficiently general,
then $W$ defines an exact symplectic Lefschetz fibration
with respect to the cylindrical K\"{a}hler form on $(\bCx)^2$,
and one can associate the directed Fukaya category
$
 \Fuk W
$
whose set of objects is a distinguished basis of vanishing cycles and
whose spaces of morphisms are Lagrangian intersection Floer complexes
\cite{Seidel_VC2, Seidel_PL}.
 \item
Let $X$ be the two-dimensional toric Fano stack
associated with the stacky fan
whose one-dimensional cones are generated by vertices of $\Delta$.
The derived category $D^b \coh X$
of coherent sheaves on $X$ has a full strong exceptional collection $(E_i)_i$
consisting of line bundles \cite{Borisov-Hua},
which induces a derived equivalence
with the category of finitely-generated modules
over the total morphism algebra $\bC \Gamma = \End(\bigoplus_i E_i)$.
The full subcategory of the enhanced derived category
of $\bC \Gamma$ consisting of simple modules
will be denoted by $\scC$.
 \item
Let $(G, D)$ be a pair of a consistent dimer model
and a perfect matching on it,
whose characteristic polygon coincides with $\Delta$.
One can associate a directed $A_\infty$-category $\scA$
with such a pair,
and there is a quasi-equivalence
$$
 \scA \cong \scC
$$
of $A_\infty$-categories
for a suitable choice of an exceptional collection
on $X$
\cite{Ishii-Ueda_DMEC, Futaki-Ueda_A-infinity}.
\end{itemize}
A dimer model is a bicolored graph on a real 2-torus
which encodes the information of a quiver with potential.
See e.g. \cite{Futaki-Ueda_A-infinity} and references therein
for basic definitions on dimer models.
In this paper,
we deal only with hexagonal dimer models
appearing in \cite{Ueda-Yamazaki_NBTMQ}.

The following conjecture is motivated by
\cite{Feng-He-Kennaway-Vafa}:

\begin{conjecture}[{\cite[Conjecture 6.2]{Ueda-Yamazaki_toricdP}}]
 \label{conj:FukBT}
Let $\Delta$ be a lattice polygon
containing the origin in its interior.
Then for a suitable choice of
\begin{itemize}
 \item
a Laurent polynomial $W$
whose Newton polygon coincides with $\Delta$, and
 \item
a distinguished basis of vanishing cycles on $W^{-1}(0)$,
\end{itemize}
there is a bicolored graph $Y$ on $W^{-1}(0)$
such that
\begin{itemize}
 \item
an edge of $Y$ corresponds to an intersection of vanishing cycles,
 \item
a node of $Y$ corresponds to a holomorphic disk
bounded by vanishing cycles,
 \item
the color of a node corresponds to the sign of the $A_\infty$-operation
determined by the disk,
 \item
the image of $Y$ by the argument map
$$
\begin{array}{cccc}
 \Arg : & (\bCx)^2 & \to & \bR^2 / \bZ^2 \\
 & \vin & & \vin \\
 & (x, y) & \mapsto & \dfrac{1}{2 \pi}(\arg x, \arg y)
\end{array}
$$
is a consistent dimer model $G$,
 \item
the order on the distinguished basis of vanishing cycle
gives a perfect matching $D$, and
 \item
the characteristic polygon of the pair $(G, D)$
coincides with $\Delta$.
\end{itemize}
\end{conjecture}

We prove the following in this paper:

\begin{theorem} \label{th:main}
Conjecture \ref{conj:FukBT} holds
if $\Delta$ is a triangle.
\end{theorem}

As a corollary,
one obtains a torus-equivariant version
\begin{equation} \label{eq:equiv_hms}
 D^b \coh^\bT X \cong D^b \Fuk \Wtilde
\end{equation}
of homological mirror symmetry
\cite{Kontsevich_HAMS, Kontsevich_ENS98}
for two-dimensional toric Fano stacks
of Picard number one.
Here
$
 D^b \coh^\bT X
$
is the equivariant derived category of coherent sheaves on $X$
with respect to the algebraic torus $\bT$ acting on $X$,
and $\Wtilde$
is the pull back of $W$ to the universal cover
of the torus.
The non-equivariant version of \eqref{eq:equiv_hms}
for weighted projective planes
is due to \cite{Seidel_VC2, Auroux-Katzarkov-Orlov_WPP}.

The organization of this paper is as follows:
In Section \ref{sc:wpp},
we recall the construction of two-dimensional toric Fano stacks
from lattice triangles and discuss its relation with
weighted projective planes.
In Section \ref{sc:vc},
we describe vanishing cycles of $W$ following
\cite{Auroux-Katzarkov-Orlov_WPP}
closely.
In Section \ref{sc:coamoeba},
we study the behavior of vanishing cycles
under the argument map and prove Theorem \ref{th:main}.

\section{Triangles and weighted projective planes}
 \label{sc:wpp}


Let
$$
 \Delta = \Conv \{ v_1, v_2, v_3 \}
$$
be a lattice triangle in $\bR^2$
containing the origin in its interior.
The toric Fano stack $X$ associated with $\Delta$ is
defined as the quotient stack
$$
 X = [ (\bC^3 \setminus 0) / K]
$$
of the complement of the origin in $\bC^3$
by the natural action of the kernel
$$
 K = \Ker(\phi \otimes \bCx : (\bCx)^3 \to (\bCx)^2 ),
$$
where
$$
\begin{array}{cccc}
 \phi : & \bZ^3 & \to & \bZ^2 \\
 & \vin & & \vin \\
 & e_i & \mapsto & v_i
\end{array}
$$
is the homomorphism of abelian groups
sending the $i$-th coordinate vector $e_i$
to the vertex $v_i$ of $\Delta$
for $i = 1, 2, 3$.

If $\phi$ is surjective,
then $K$ is isomorphic to $\bCx$ and
the action of $\alpha \in \bCx \cong K$
on $\bC^3$ is given by
$(x, y, z) \mapsto (\alpha^a x, \alpha^b y, \alpha^c z)$
for some relatively prime positive integers
$a$, $b$, and $c$.
The resulting stack $[(\bC^3 \setminus 0) / K]$
is the weighted projective plane
$$
 X = \bP(a, b, c),
$$
and any weighted projective plane
with relatively prime weights can be obtained in this way
by setting $\phi$ to be the natural projection
$$
 \phi : \bZ^3 \to \coker(\varphi) \cong \bZ^2
$$
to the cokernel of
$$
\begin{array}{cccc}
 \varphi : & \bZ & \to & \bZ^3 \\
 & \vin & & \vin \\
 & 1 & \mapsto & (a, b, c).
\end{array}
$$
If $d = \gcd(a, b, c) \ne 1$,
then the derived category of coherent sheaves on $\bP(a, b, c)$
is a direct sum
$$
 D^b \coh \bP(\bsa, \bsb, \bsc)
  \cong (D^b \coh \bP(\bsa', \bsb', \bsc'))^{\oplus \bsd},
 \qquad (\bsa, \bsb, \bsc) = (\bsa' \bsd, \bsb' \bsd, \bsc' \bsd),
$$
and the mirror of $\bP(a, b, c)$ is a disjoint union
of $\bsd$ copies of the mirror for $\bP(\bsa', \bsb', \bsc')$.

If the map $\phi$ is not surjective,
then one can factor $\phi$ as
$$
 \phi = \phi_2 \circ \phi_1 : \bZ^3 \xto{\phi_1} \bZ^2 \xto{\phi_2} \bZ^2
$$
where $\phi_1$ is the surjection to $\Image(\phi) \cong \bZ^2$
and $\phi_2$ is the inclusion of $\Image(\phi)$ to $\bZ^2$.
One obtains an exact sequence
$$
 1 \to K_1 \to K \to K_2 \to 1,
$$
where
$
 K_1 = \Ker (\phi_1 \otimes \bCx)
$
and
$
K_2 = \Ker (\phi_2 \otimes \bCx),
$
and $X = [ (\bC^3 \setminus 0) / K]$ is the quotient stack
$$
 X = [\bP(a, b, c) / K_2]
$$
for the weight $(a, b, c)$ such that
$
 \bP(a, b, c) = [(\bC^3 \setminus 0) / K_1].
$

\section{Vanishing cycles for triangles}
 \label{sc:vc}

\begin{figure}
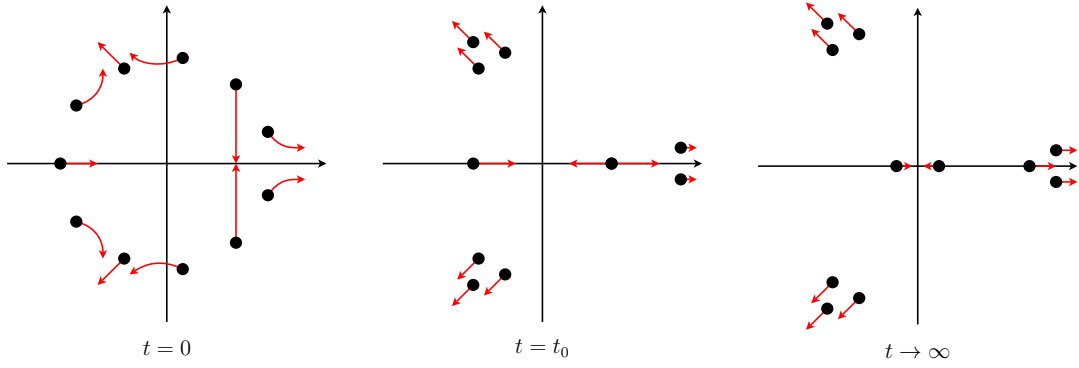

\begin{minipage}{.3 \linewidth}
\centering
\scalebox{.7}{\input{x-plane_vc3.pst}}
\end{minipage}
\begin{minipage}{.3 \linewidth}
\centering
\scalebox{.7}{\input{x-plane_vc4.pst}}
\end{minipage}
\begin{minipage}{.3 \linewidth}
\centering
\scalebox{.7}{\input{x-plane_vc5.pst}}
\end{minipage}
\caption{Branch points of $\varpi_t$ for $(a, g, h) = (3, 3, 2)$}
\label{fg:x-plane_vc}
\end{figure}

Let $\Delta \subset \bR^2$ be a lattice triangle
which contains the origin in its interior.
One can choose an $\SL_2(\bZ)$-transformation
to set
$$
 \Delta = \Conv \{ (a, 0), (b, c), (-d, -e) \}
$$
where $a, c, d, e$ are positive and
$b$ is non-negative.
Let
$$
 W = x^a + x^b y^c + \frac{1}{x^d y^e}
$$
be a Laurent polynomial
whose Newton polygon coincides with $\Delta$ and
consider the diagram
\vspace{7mm}
$$
\begin{psmatrix}[colsep=1.5]
 \bCx & \bC \times \bCx & \bC
\end{psmatrix}
\psset{shortput=nab,arrows=->,labelsep=3pt}
\ncline{1,1}{1,2}_{\varpi}
\ncline{1,2}{1,3}_{\psi}
\ncarc[arcangle=30]{1,1}{1,3}^{W = \psi \circ \varpi}
$$
where
$$
 \varpi(x, y) = (W(x, y), x)
$$
and
$$
 \psi(w, x) = w.
$$
For general $t \in \bC$, the map
$$
 \scE_t \xto{\varpi_t} \scS_t
$$
from $\scE_t = W^{-1}(t)$ to $\scS_t = \psi^{-1}(t)$
is a $(c + e)$-fold cover of $\scS_t$,
which can naturally be identified with the $x$-plane.
The fiber of $\varpi_t$ is defined by
$$
 x^{b+d} y^{c+e} + (x^a - t) x^d y^e + 1 = 0,
$$
which can be written as
$$
 A y^{c+e} + B y^e + 1 = 0
$$
where
$A = x^{b+d}$ and $B = (x^a - t) x^d$.
The critical points of $\varpi_t$ are defined by
\begin{align} \label{eq:crit_varpi}
\left\{
\begin{aligned}
 A y^{c+e} + B y^e + 1 &= 0, \\
 (c+e) A y^{c+e-1} + e B y^{e-1} &= 0.
\end{aligned}
\right.
\end{align}
By eliminating $y$ from \eqref{eq:crit_varpi},
one obtains
$$
 (-1)^{g} \frac{c^c e^e}{g^{g}} (x^a - t)^g x^h = 1
$$
as the defining equation for the critical values of $\varpi_t$
where $g = c + e$ and $h = c d - b e$.
The set
$$
 D(t) = \lc x \in \bCx \, \left| \,
  (-1)^{g} \frac{c^c e^e}{g^{g}} (x^a - t)^g x^h = 1 \right. \rc
$$
consists of $a g + h$ points
for general $t$,
which becomes singular when
\begin{align*}
 x^a &= \frac{h}{a g + h} t
\end{align*}
and
\begin{equation} \label{eq:critv_W}
 t^{g + \frac{h}{a}} = \frac{1}{c^c e^e}
  \lb 1 + \frac{h}{a} \rb^g \lb 1 + \frac{a g}{h} \rb^{\frac{h}{a}}.
\end{equation}
Assume $k := \gcd(a, h) = 1$,
so that \eqref{eq:critv_W} has $a g + h$ solutions,
which is equal to the area $|\Delta|$ of $\Delta$.
The $k \ne 1$ case can be reduced
to this case by the $k$-fold cover $x \mapsto x^k$
of the $x$-plane.
The set of solutions of \eqref{eq:critv_W} is
the set of critical values of $W$.
We choose the straight line segments
from the origin to the critical values of $W$
as a distinguished set of vanishing paths.
The corresponding vanishing cycles can be computed
by studying the behavior of the branch points of $\varpi_t$
along the vanishing paths.
Let $t_0$ be the unique positive real critical value of $W$ and
consider the behavior of the branch points of $\varpi_t$
as one varies $t$ from zero to infinity along the positive real axis.
Branch points at $t = 0$ are distributed
on a circle centered at the origin,
and their arguments are given by $\frac{(g + 2 n) \pi}{a g + h}$
for $n = 1, \ldots, a g + h$.
As $t$ goes from zero to $t_0$,
the branch points with arguments $\pm \frac{g \pi}{a g + h}$
come close to each other and merge on the real line.
As $t$ goes from $t_0$ to infinity,
the merged branch points split into two again,
and the whole set of branch points are divided into $a + 1$ groups;
one group consists of $h$ branch points coming close to the origin,
and each of the remaining $a$ groups
consists of $g$ branch points going off to infinity.
Figure \ref{fg:x-plane_vc} shows this behavior
for $(a, g, h) = (3, 3, 2)$.
It follows that
the vanishing cycle $C_0 \subset W^{-1}(0)$ of $W$
along the straight line segment from the origin to $t_0$
lies above the matching path
obtained as the trajectory of two branch points of $\varpi_t$
whose arguments are $\pm \frac{g \pi}{a g + h}$ at $t = 0$.

Let $\phi_0 : \bZ^2 \to \bZ^2$ be
the linear map represented by the matrix
$$
\begin{pmatrix}
 a + d & b + d \\
 e & c + e
\end{pmatrix},
$$
and
$$
 K_0 = \Ker(\phi_0 \otimes \bCx)
$$
be the kernel of the homomorphism
$$
 \phi_0 \otimes \bCx : (\bCx)^2 \to (\bCx)^2.
$$
Then an element $(\alpha, \beta) \in K_0$ gives a map
$$
\begin{array}{ccc}
 W^{-1}(t) & \to & W^{-1}(\alpha^d \beta^e t) \\
 \vin & & \vin \\
 (x, y) & \mapsto & (\alpha x, \beta y),
\end{array}
$$
which induces a free transitive action of $K_0$
on the distinguished basis of vanishing cycles on $W^{-1}(0)$
along straight line segments
from the origin to critical values of $W$.

\section{Coamoebas and vanishing cycles}
 \label{sc:coamoeba}

Recall from \cite[Theorem 7.1]{Ueda-Yamazaki_NBTMQ}
that the coamoeba of
$$
 W^{-1}(0) = \{ (x, y) \in (\bCx)^2 \mid
  1 + x^{a+d} y^e + x^{b+d} y^{c+e} = 0 \}
$$
is given by the pull-back of the coamoeba of
$$
 \{ (x, y) \in (\bCx)^2 \mid 1 + x + y = 0 \}
$$
shown in Figure \ref{fg:Z1_coamoeba}
by the map
$
 \psi \otimes (\bR / \bZ) : \bR^2 / \bZ^2 \to \bR^2 / \bZ^2,
$
where $\psi : \bZ^2 \to \bZ^2$ is the linear map
represented by the matrix
$$
\begin{pmatrix}
 p & q \\
 r & s
\end{pmatrix}
=
\begin{pmatrix}
 a + d & e \\
 b + d & c + e
\end{pmatrix}.
$$
Figure \ref{fg:coamoeba_1} shows a part of
the coamoeba,
and the entire coamoeba is obtained
by gluing $|\Delta| = p s - r q$ copies of it.

\begin{figure}
\begin{minipage}{.5 \linewidth}
\centering
\input{Z1_coamoeba.pst}
\caption{The coamoeba of $x + y + 1 = 0$}
\label{fg:Z1_coamoeba}
\end{minipage}
\begin{minipage}{.5 \linewidth}
\centering
\scalebox{.9}{\input{coamoeba_1.pst}}
\caption{A part of the coamoeba of $W^{-1}(0)$}
\label{fg:coamoeba_1}
\end{minipage}
\end{figure}

The coamoeba is the union of open triangles and their vertices,
and the inverse image of the set of vertices
of the coamoeba of $1 + x + y = 0$
is the real part of $1 + x + y = 0$.
It is parametrized as
$$
\left\{
\begin{aligned}
 x &= t, \\
 y &= - t - 1,
\end{aligned}
\right.
$$
and divided into three parts
$$
 t < -1, \quad
 -1 < t < 0, \quad \text{and} \quad 
 t > 0.
$$
It follows that the inverse images of vertices
of the coamoeba of $W^{-1}(0)$ is parametrized as
$$
 (\psi \otimes \bCx)^{-1}(t, -1-t).
$$
Since $\psi^{-1}$ is given by the matrix
$$
\frac{1}{|\Delta|}
\begin{pmatrix}
 c + e & - e \\
 - b - d & a + d
\end{pmatrix},
$$
the $x$-projection
of the inverse images of vertices of the coamoeba
is parametrized as
\begin{align*}
 x^{|\Delta|} &= \frac{t^{c + e}}{(- 1 - t)^e}.
\end{align*}
By studying the behavior of the function
$$
 f(t) = \frac{t^{c+e}}{(-1-t)^{e}},
$$
one can see that the $x$-projections
of the inverse images of vertices of the coamoeba
corresponding to the vertex $(\frac{1}{2}, 0)$
of the coamoeba of $x + y + 1$
are half lines
from the branch points of $\varpi_0$ to infinity
with constant arguments.
The $x$-projections of inverse images of other vertices
of the coamoeba are half lines
from the origin to infinity with constant arguments.
The fiber $W^{-1}(0)$ is obtained
by gluing $c+e$ copies of the $x$-plane
which are cut into $2 |\Delta|$ pieces
along these half lines.

\begin{figure}
\begin{minipage}{.5 \linewidth}
\centering
\input{x-plane_vc1.pst}
\caption{Three triangles and a part of the vanishing cycle $C_0$}
\label{fg:x-plane_vc1}
\end{minipage}
\begin{minipage}{.5 \linewidth}
\centering
\input{x-plane_vc2.pst}
\caption{Three other triangles and
the other part of the vanishing cycle $C_0$}
\label{fg:x-plane_vc2}
\end{minipage}
\begin{minipage}{.5 \linewidth}
\centering
\input{coamoeba_2.pst}
\caption{The vanishing cycle $C_0$ on the coamoeba}
\label{fg:coamoeba_2}
\end{minipage}
\begin{minipage}{.5 \linewidth}
\centering
\input{coamoeba_3.pst}
\caption{The dimer model}
\label{fg:coamoeba_3}
\end{minipage}
\end{figure}

Now consider six triangles in Figure \ref{fg:coamoeba_2}
which are adjacent to two triangles in Figure \ref{fg:coamoeba_1}.
The corresponding pieces of the copies of the $x$-plane
are shown in Figure \ref{fg:x-plane_vc1} and Figure \ref{fg:x-plane_vc2}.
The discussion in Section \ref{sc:vc} shows that
the vanishing cycle $C_0 \subset W^{-1}(0)$
along the straight line segment
from the origin to the positive real critical value
is obtained by gluing the curves
in Figure \ref{fg:x-plane_vc1} and Figure \ref{fg:x-plane_vc2}
connecting branch points with arguments $\pm \frac{s}{\Delta} \pi$.
The argument projection of $C_0$ is
shown in Figure \ref{fg:coamoeba_2},
which naturally corresponds
to a face of the hexagonal dimer model $G$
shown in Figure \ref{fg:coamoeba_3}.
Other vanishing cycles are obtained from $C_0$
by the action of $K_0$
as described in Section \ref{sc:vc},
so that the argument projection induces
a natural bijection
between the distinguished basis of vanishing cycles of $W$
and the set of faces of $G$.
One can easily see that an edge of $G$ corresponds to
an intersection of vanishing cycles under this bijection,
and a node of $G$ gives a holomorphic triangle
which contributes to the $A_\infty$-operation $\frakm_2$
on the Fukaya category.
When $X$ is a weighted projective plane, a comparison
with the discussion in \cite{Auroux-Katzarkov-Orlov_WPP} shows that
the color of the node matches the sign in the $A_\infty$-operation,
and the ordering on the distinguished basis of vanishing cycles
defines an internal perfect matching of $G$.
The toric Fano stack associated with a general lattice triangle
can be obtained from the weighted projective plane
as a toric orbifold, and
Theorem \ref{th:main} is proved.

{\em Acknowledgment}:
M.~F. is supported by Grant-in-Aid for Young Scientists (No.19.8083).
K.~U. is supported by Grant-in-Aid for Young Scientists (No.18840029).

\bibliographystyle{amsalpha}
\bibliography{bibs.bib}

\noindent
Masahiro Futaki

Graduate School of Mathematical Sciences,
The University of Tokyo,
3-8-1 Komaba Meguro-ku Tokyo 153-8914, Japan

{\em e-mail address}\ : \  futaki@ms.u-tokyo.ac.jp

\ \\

\noindent
Kazushi Ueda

Department of Mathematics,
Graduate School of Science,
Osaka University,
Machikaneyama 1-1,
Toyonaka,
Osaka,
560-0043,
Japan.

{\em e-mail address}\ : \  kazushi@math.sci.osaka-u.ac.jp
\ \vspace{0mm} \\

\end{document}